\documentclass[11pt]{amsart}    

\usepackage{graphicx}
\usepackage{color}
\definecolor{Myred}{cmyk}{0.0,1.0,1.0,0.00}
\definecolor{Mygreen}{rgb}{0,0.4,0}
\usepackage{enumerate}

\usepackage{amsfonts,amsmath,latexsym,amssymb,amsthm}
\usepackage{hyperref}
\usepackage{geometry}
\usepackage{txfonts}

\newtheorem{theorem}{Theorem}

\theoremstyle{remark}
\newtheorem{remark}{Remark}

\begin{document}

\title{Lieb-Thirring type estimates for Dirichlet Laplacians on spiral-shaped domains}

\author{Juan Bory-Reyes}
\address{ESIME-Zacatenco, Instituto Polit\'{e}cnico Nacional, M\'{e}xico, CDMX. 07738. M\'{e}xico}
\email{juanboryreyes@yahoo.com}

\author{Diana Barseghyan (Schneiderov\'{a})}  
\address{Department of Mathematics, University of Ostrava,  30. dubna 22, 70103 Ostrava, Czech Republic}
\email{diana.schneiderova@osu.cz}            

\author{Baruch Schneider}
\address{Department of Mathematics, University of Ostrava,  30. dubna 22, 70103 Ostrava, Czech Republic}
\email{baruch.schneider@osu.cz} 


\maketitle

\begin{abstract}
In this paper we derive Lieb-Thirring estimates for eigenvalues of Dirichlet Laplacians below the threshold of the essential spectrum on asymptotically Archimedean spiral-shaped regions.
\end{abstract}

\section{Introduction} \label{s: intro}

The dynamics of quantum particles confined in unbounded regions of various shapes is not only important from a physical point of view, but also from a mathematical point of view, where important connections between spectral properties of the corresponding Hamiltonians and the confinement geometry are described. We refer the reader to \cite{EK15} for background and further references. 
While some geometric perturbations types, such as bends or twists of straight tubes, both local and periodic, have been
deeply studied mathematically, some others have escaped attention. This is especially the case for spiral structures, which appear in physics, for example, as waveguides for cold atoms \cite{JLXZY15}. A mathematical analysis of Dirichlet Laplacians in spiral-shaped regions has recently been presented in \cite{ET21} and \cite{BE23}. There are many different types of spiral regions. A decisive factor for the spectral properties is the behaviour of the coil width as we follow the spiral from the center to infinity; among those for which this quantity is monotonous, we can distinguish spirals expanding, asymptotically Archimedean, and shrinking -- definitions will be given below. The extreme cases are when the coil width diverges or shrinks to zero, then the spectrum covers the whole positive halfline or it is purely discrete, respectively. The situation is more interesting when the region is asymptotically Archimedean in the sense that the coil width has a finite and nonzero limit.

In the present paper we consider the asymptotically Archimedean spiral-shaped regions for which the Dirichlet Laplacian spectrum consists of the essential part and the eigenvalues below the threshold of the essential spectrum.

Our purpose is to obtain the bounds for the moments of these eigenvalues in terms of the geometric properties of the region in the spirit of the estimates derived by Lieb and Thirring, Berezin, Lieb, and Li and Yau, cf.~\cite{LT76, Be72a, Be72b, Li73, LY83} and the monograph \cite{FLW22}.

The cited results concern situations in which the motion is restricted to a finite domain or is governed by a Schr\"odinger operator with a finite classically allowed domain. However, it has long been known that the spectrum of the Dirichlet Laplacian can be discrete, in whole or in part, even when the phase space volume involved is infinite with respect to regions of tube or cusp-type shapes; the corresponding nonclassical eigenvalue asymptotic behaviour has been studied, e.g. in \cite{Re48, Ro72, Ro73, Si83, vdB84, vdB92} with the focus on the spectral counting function. Other moments of the eigenvalue distribution have also been considered, e.g. in \cite{ELW04, GW11, BE13, EB14, BK19}.

The analysis of spiral-shaped regions that we will present is close to the above mentioned case. Depending on the function $d(\cdot)$ describing the width of the spiral, the spectrum can be purely discrete, purely essential, or the combination of the essential spectrum with the eigenvalues below its threshold. As mentioned above, we consider the situation of the asymptotically Archimedean spiral-shaped regions.

Our main result, an estimate of the eigenvalue moments below the threshold of the essential spectrum in terms of the geometric properties of the spiral, is presented and proved as Theorem~\ref{Main} in Sec.~\ref{Main result}. 
Sec.~\ref{Absence of the infinite discrete spectrum} presents some examples of the spirals for which the discrete spectrum is absent. Finally, Sec.~\ref{Local perturbations of the Archimedean spiral} deals with the eigenvalue bounds in the case when the curve is the local perturbation of the Archimedean spiral.

\section{Preliminaries} \label{s: prelim}

First, let us describe the geometry of spiral-shaped regions. Let $\mathcal{C}$ be the graph of an increasing function $r : \mathbb{R}_+ \to \mathbb{R}_+$ with $r(0)=0$ and $\lim_{\theta\to\infty}r(\theta)=\infty$, that is, the family of points $(r(\theta),\theta)_{\theta\in(0, \infty)}$ in the polar coordinates.

The assumed monotonicity of $r$ means that $\mathcal{C}$ does not intersect itself which means, in particular, that the width function
\begin{equation}\label{a}
a(\theta):=\frac{1}{2\pi}\big(r(\theta+2\pi)-r(\theta)\big)
\end{equation}
is positive for any $\theta\ge 0$. A spiral curve $\mathcal{C}$ is called \emph{simple} if the corresponding $a(\cdot)$ is monotonous, and \emph{expanding} or \emph{shrinking} if this function is increasing or decreasing, respectively, in $\theta$ away from a neighborhood of the origin; 
these qualifications are labeled as \emph{strict} if $\lim_{\theta\to\infty}a(\theta)=\infty$ and $\lim_{\theta\to\infty}a(\theta)=0$, respectively. A simple spiral lying between these two extremes, for which the limit is finite and nonzero, is called \emph{asymptotically Archimedean}.
The main object of our interest is the two-dimensional Laplace operator $H_\Omega$ with the Dirichlet condition imposed at the boundary represented by the curve $\mathcal{C}$, defined in the standard way \cite[Sec.~XIII.15]{RS78} on the open set $\Omega=\Omega_\mathcal{C}=\mathbb{R}^2\setminus\mathcal{C}$. As shown in \cite{ET21}, the spectral properties of such operators on simple spiral regions depend strongly on the function $a(\cdot)$. For strictly expanding regions the spectrum is purely essential and covers the halfline $\mathbb{R}_+$. On the other hand, if the spiral 
$\mathcal{C}$ is strictly shrinking, the spectrum of $H_\Omega$ is purely discrete. 

In the intermediate case of asymptotically Archimedean spirals the spectrum may be more complicated. Its essential part covers the interval 
$\big[\frac14(\lim_{\theta\to\infty}a(\theta))^{-2}, \infty\big)$. The discrete part can be empty, as in the case of the pure Archimedean spiral, if $a(\theta)\equiv a_0$, but also infinite, accumulating at the threshold of $\sigma_\mathrm{ess}(H_\Omega)$ if the spiral is the the corresponding perturbation of the pure Archimedean spiral \cite{ET21}.

In this paper we consider the case of asymptotically Archimedean spirals and obtain bounds for the eigenvalues below the threshold of the essential spectrum. Using the dimension reduction technique of Laptev and Weidl \cite{LW00}, we derive Lieb-Thirring-type inequalities for eigenvalue moments in terms of the curvature and the width function of the curve $\mathcal{C}$.

A useful way to characterize the region $\Omega$, possibly with the exception of a neighborhood of the origin of the coordinates, is to use the Fermi (or parallel) coordinates, that is, a locally orthogonal system in which the Cartesian coordinates of $\mathcal{C}$ are written as
\begin{eqnarray}
\nonumber x_1(\theta, u) = r(\theta)\cos\theta-\frac{u}{\sqrt{\dot{r}(\theta)^2+r(\theta)^2}}
(\dot{r}(\theta) \sin\theta+r(\theta)\cos\theta),\\ \label{Fermi}
x_2(\theta, u)=r(\theta)\sin\theta+ \frac{u}{\sqrt{\dot{r}(\theta)^2+r(\theta)^2}}(\dot{r}(\theta)\cos\theta-r(\theta)\sin\theta),
\end{eqnarray}
where $u$ measures the distance of $(x_1, x_2)$ from $\mathcal{C}$. A natural counterpart of the variable $u$ is the arc length of the spiral given by
\begin{equation}\label{s}
s(\theta)=\int_0^\theta\sqrt{\dot{r}(t)^2+r(t)^2}\,d t.
\end{equation}
We use relations \eqref{Fermi} and \eqref{s} to parameterize the region $\Omega$ with the coordinates $(s,u)$, possibly with the exception of a finite central part, like
\begin{equation}\label{parametrization}
\Omega_1=\Omega\cap\{(s, u): s>s_0\}= \big\{x(s, u): s >s_0,\: u \in(0, d(s))\big\}\,;
\end{equation}
where $s_0>0$ is a number depending on the curve $\mathcal{C}$ characterizing the excluded part, and $d(s)$ is the length of the inward normal starting from the point $x(s, 0)$ of $\mathcal{C}$ to the intersection with the previous coil of the spiral.

One more quantity associated with the spiral that we will need to state the result is its \emph{curvature} which is in terms of the angular variable given by
\begin{equation}\label{curvature}
\gamma(\theta)= \frac{r(\theta)^2+ 2\dot{r}(\theta)^2- r(\theta)\ddot{r}(\theta)}{
(r(\theta)^2+\dot{r}(\theta)^2)^{3/2}},
\end{equation}
assuming, of course, that the derivatives make sense. Using the pull-back, $s\mapsto\theta(s)$, of the map \eqref{s}, we can alternatively express it as a function of the arc length $s$. We abbreviate $\gamma(\theta(s))$ to $\gamma(s)$.

\section{Main result}\label{Main result}

To be more specific, we assume that function $r$, in addition to its monotonicity and the condition $r(0)=0$, is subject to the
condition:
\begin{enumerate}[(a)]
\setlength{\itemsep}{0pt}\item
$r(\theta)= a_0 \theta+ \rho(\theta),\,a_0>0;$ \,\,
$\lim_{\theta\to\infty}\rho(\theta)=0$.  \label{assa}
 \end{enumerate}
Under this assumption it is easy to see that $\mathcal{C}$ is asymptotically Archimedean.

The parallel coordinates can be used in regions generated by asymptotic Archimedean spirals, provided we exclude a suitable central region. This could be checked by the facts that $\gamma(s)=\gamma(s(\theta)) = \mathcal{O}((s(\theta))^{-1})$ as $\theta\to\infty$ vanishes to zero, given by equation (\ref{curvature}), and $d(s)$ is upper bounded by $2\pi a_0+\|\rho(\tau+2\pi)- \rho(\tau)\|_{L^\infty(0, \infty)}$, there exists $s_0>0$ such that $d(s)\gamma(s)<1$ holds for all $s\ge s_0$, which means, in particular, that in the corresponding part of $\Omega$ the Fermi coordinates are well defined.

It is known that the essential spectrum of $H_\Omega$ coincides with the semi-interval $\left[\frac{1}{4a_0^2}, \infty\right)$ \cite{ET21}. 
The problem we are interested in is the bounds for the eigenvalues of the operator $H_\Omega$ below the $\frac{1}{4a_0^2}$, more precisely the moments of the negative eigenvalues of the operator $H_\Omega- \frac{1}{4a_0^2}$. 
\begin{theorem}\label{Main}
Let $\Omega=\Omega_\mathcal{C}$ be a domain determined by an asymptotically Archimedean spiral $\mathcal{C}$ satisfying the assumption \eqref{assa}, and let $H_\Omega$ be the corresponding Dirichlet Laplacian. Then the following inequality holds for any $\sigma\ge\frac{3}{2}$,
\begin{align}\nonumber
&\mathrm{tr}\left(H_\Omega-\Lambda\right)_-^\sigma \\
\label{theorem}&\le
\,\frac{2L_{\sigma,1}^{\mathrm{cl}}}{\pi}\int_{s_0}^\infty \sqrt{\widetilde{W}(s)+\frac{1}{4a_0^2}}\, \left(\widetilde{W}(s)+\frac{1}{4a_0^2}-\left(\frac{\pi }{d(s)} \right)^2 \right)_+^{\sigma+\frac{1}{2}} d(s)\,d s
 + \frac{L_{\sigma, 2}^{\mathrm{cl}}\,\mathrm{vol}(\Omega_2)}{4^{\sigma+1/2} a_0^{2\sigma+2}},
\end{align}
where $\Omega_2$ is the finite subset of $\Omega$  with respect to the arc length $s\le s_0$, and $L_{\sigma,1}^{\mathrm{cl}}$ is the semiclassical constant,
\begin{equation} \label{LTconstant}
L_{\sigma,1}^{\mathrm{cl}} := \frac{\Gamma(\sigma+1)}{\sqrt{4\pi} \Gamma(\sigma+\frac32)}\,;
 \end{equation}
the function $\widetilde{W}$ in \eqref{theorem} is given by
\begin{equation}\label{W}
\widetilde{W}(s):=\frac{\gamma^2(s)}{4(1-\gamma(s)d(s))^2}+\frac{d(s)|\ddot\gamma(s)|}{2(1-\gamma(s)d(s))^3}+\frac{5}{4}
\frac{d(s)^2|\dot\gamma(s)|^2}{(1-d(s)\gamma(s))^4}.
\end{equation}
\end{theorem}
\medskip

\begin{proof}[Proof of Theorem \ref{Main}]
As mentioned above, the assumption \eqref{assa} allows us to use the Fermi coordinate parametrization (\ref{parametrization}) of $\Omega$ with $s_0>0$. We use the Neumann bracketing method \cite[Sec.~XIII.15]{RS78}, which yields
\begin{equation}\label{bracketing}
H_\Omega-\frac{1}{4a_0^2}\ge \left(H_{\Omega_1}-\frac{1}{4a_0^2}\right)\oplus \left(H_{\Omega_2}-\frac{1}{4a_0^2}\right),
\end{equation}
where $H_{\Omega_1}$ and $H_{\Omega_2}$ are the restrictions of $H_\Omega$ with respect to the regions $\Omega_1\subset\Omega$ corresponding to the arc lengths $s>s_0$ and $\Omega_2:=\Omega\setminus\bar{\Omega}_1$, both having the additional Neumann condition at $s=s_0$.

Consider first the operator $H_{\Omega_1}$. According to \cite{EK15}, the coordinates (\ref{Fermi}) allow us to pass from $H_{\Omega_1}$ to a unitarily equivalent operator $\widetilde{H}_{\widetilde{\Omega}_1}$ defined on the region $\widetilde{\Omega}_1:=\{(s, u): s>s_0,\, 0<u<d(s)\}$ with the Neumann boundary condition at $s=s_0$ and the Dirichlet condition on the rest of the boundary of $\Omega_1$ as follows,
\begin{equation}\label{unitary eq.}
\big(\widetilde{H}_{\widetilde{\Omega}_1}\psi\big)(s,u) = -\left(\frac{\partial}{\partial s}\frac{1}{(1-u\gamma(s))^2}\frac{\partial\psi}{\partial s}\right)(s,u)-\frac{\partial^2\psi}{\partial u^2}(s,u)+W(s, u)\psi(s,u)\,,
\end{equation}
where
$$
W(s, u):=\frac{\gamma^2(s)}{4(1-u\gamma(s))^2}+\frac{u\ddot\gamma(s)}{2(1-u\gamma(s))^3}+\frac{5}{4}
\frac{u^2\dot\gamma(s)^2}{(1-u\gamma(s))^4}.
$$
It is straightforward to check that
\begin{equation}\label{est.}
\widetilde{H}_{\widetilde{\Omega}_1}\ge H_0,
\end{equation}
where the operator $H_0=-\Delta-\widetilde{W}$ acts on $L^2(\widetilde{\Omega}_1)$ and satisfies the same boundary conditions as $\widetilde{H}_{\widetilde{\Omega}_1}$, and $\widetilde{W}$ is given by (\ref{W}).

So we can get the desired result by setting a bound on the trace of the operator $H_0$.

An original impulse for this investigation came from \cite{LW00} and the use of a variational argument to reduce the problem to a Lieb-Thirring inequality with an operator-valued potential.
Given a function $g\in\,C^\infty (\widetilde{\Omega}_1)$ with a zero trace on the set $\{s\in(s_0, \infty),\, u=0,\,d(s)\}$, we have
\begin{align*}
& \int_{\widetilde{\Omega}_1}\left(|\nabla\,g(s,u)|^2-\left(\widetilde{W}(s)+\frac{1}{4a_0^2}\right)|g(s,u)|^2\right)\,d s\,d u \\ & =\int_{\widetilde{\Omega}_1}\Big|\frac{\partial g}{\partial s}(s,u)\Big|^2\,d s\,d u
+\int_{s_0}^\infty\,\mathrm{d}s \int_0^{d(s)}\bigg(\Big|\frac{\partial g}{\partial u}(s,u)\Big|^2
- \left(\widetilde{W}(s)+\frac{1}{4a_0^2}\right) |g(s,u)|^2\bigg)\,d u \\
& =\int_{\widetilde{\Omega}_1}\Big|\frac{\partial g}{\partial s}(s,u)\Big|^2\,d s\,d u +\int_{s_0}^\infty\left\langle\,H(s,\widetilde{W}(s))g(s,\cdot), g(s,\cdot)\right\rangle_{L^2(0, d(s))}
\,d s\,,
\end{align*}
where $H(s,\widetilde{W}(s))$ is the one-dimensional operator
\begin{equation}\label{HW}
H(s, \widetilde{W}(s))=-\frac{\mathrm{d}^2}{\mathrm{d}u^2} -\widetilde{W}(s)-\frac{1}{4a_0^2}
\end{equation}
defined on $L^2(0, d(s))$ with Dirichlet conditions at $u=0$ and $u=d(s)$.

Next we consider the complement of $\widetilde{\Omega}_1$ to the halfplane $\{s\ge s_0, u\in\mathbb{R}\}$ and denote its interior as $\widetilde{\Omega}_1^{\mathrm{c}}$. We take arbitrary functions $g\in C^\infty(\widetilde{\Omega}_1)$ and $v\in C^\infty(\widetilde{\Omega}_1^{\mathrm{c}})$, both of which have a zero trace at $\{s\in(s_0, \infty),\, u=0,\,d(s)\}$; by extending them by zero to the complements of $\widetilde{\Omega}_1$ and $\widetilde{\Omega}_1^{\mathrm{c}}$, respectively, we can consider them as functions in the whole halfplane. Similarly, we extend $H(s, \widetilde{W}(s))$ to the operator on $L^2(\mathbb{R})$, which acts on $H(s, \widetilde{W}(s))\oplus 0$ with the zero component on 
$\mathbb{R}\setminus[0, d(s)]$. For their sum, $h=g+v$, we then have
\begin{align*}
& \|\nabla\,g\|^2_{L^2(\widetilde{\Omega}_1)} +\|\nabla\,v\|^2_{L^2(\widetilde{\Omega}_2)} -\int_{\widetilde{\Omega}_1}\left(H(s, \widetilde{W})(s)+\frac{1}{4a_0^2}\right) |g(s,u)|^2\,d s\,d u
\\ & \quad
\ge\int_{\{s>s_0, u\in\mathbb{R}\}}\Big|\frac{\partial h}{\partial s}(s,u)\Big|^2\,d s\,d u +\int_{s_0}^\infty\left\langle
\,H(s, \widetilde{W}(s))\,h(s,\cdot),\,h(s,\cdot)\right\rangle_{L^2(\mathbb{R})}
\,d s\,.
\end{align*}
The left-hand side of this inequality is the quadratic form corresponding to the direct sum of the operator $H_0-\frac{1}{4a_0^2}$ and the Laplace operator defined on $\widetilde{\Omega}_1^{\mathrm{c}}$ with the Neumann boundary conditions at $s=s_0$ and Dirichlet conditions at the rest part of the boundary, while the right-hand side is the form associated with the operator
\begin{equation}\label{Q}
Q=-\frac{\partial^2}{\partial\,s^2}\otimes\,I_{L^2(\mathbb{R})} +H(s, \widetilde{W}(s))\,,
\end{equation}
the form domain is $\mathcal{H}^1\left((s_0, \infty), L^2(\mathbb{R})\right)$. Since the Laplace operator is positive, then for any $\sigma\ge0$ the minimax principle implies
\begin{equation}\label{op.val.}
\mathrm{tr}\,\left(H_0-\frac{1}{4a_0^2}\right)_-^\sigma \le \,\mathrm{tr}\Big(-\frac{\partial^2}{\partial\,s^2}\otimes \,I_{L^2(\mathbb{R})}+ H(s, \widetilde{W}(s))\Big)_-^\sigma.
\end{equation}

Let us now use the following version of the Lieb-Thirring inequality for operator-valued potentials (the proof of which is given in the Appendix):
\begin{theorem}\label{operator valued1}
Let the operator $Q$ be given by  (\ref{HW}) and (\ref{Q}). Then the following estimate holds for any $\sigma\ge 3/2$,
\begin{equation}\label{new result}
\mathrm{tr}Q_-^\sigma
\le\, 2L_{\sigma, 1}^{\mathrm{cl}}\int_0^\infty\mathrm{tr}\,(H(s, \widetilde{W}(s)))_-^{\sigma+1/2}(s)\,d s,
\end{equation}
Furthermore, for any $\sigma\ge1/2$, the inequality is similar to the one above, but with a different constant 
\begin{equation}\label{new result.}
\mathrm{tr}Q_-^\sigma
\le\, 2r(\sigma,1) L_{\sigma, 1}^{\mathrm{cl}}\int_0^\infty\mathrm{tr}\,(H(s, \widetilde{W}(s)))_-^{\sigma+1/2}\,d s,
\end{equation}
where $r(\sigma,1)\le 2$ if $\sigma<3/2$.
\end{theorem}
  
\begin{remark}\label{s0} 
{\rm It is easy to check that the Theorem \ref{operator valued1} still
holds with $L^2(\mathbb{R}_+, L^2(\mathbb{R}))$ replacing by $L^2((\alpha, \infty), L^2(\mathbb{R}))$ with some
$\alpha\in\mathbb{R}$; the difference being only the integration interval on the right-hand side of \eqref{new result} and (\ref{new result.}) to $(\alpha, \infty)$.}
\end{remark}

Now let us return to the estimate \eqref{op.val.}. Given the Theorem \ref{operator valued1}, we have
\begin{equation}\label{op.val.new}
\mathrm{tr}\,\left(-\frac{\partial^2}{\partial\,s^2}\otimes \,I_{L^2(\mathbb{R})}+ H(s, \widetilde{W}(s))\right)_-^\sigma \le \,2L_{\sigma,1}^{\mathrm{cl}} \int_{s_0}^\infty\mathrm{tr}\,H(s, \widetilde{W}(s))_-^{\sigma+1/2}\,d s.
\end{equation}

Since the eigenvalues of $H(s,0)$ coincide with $\left(\frac{\pi j}{d(s)}\right)^2,\: j=1,2,\dots\,,$ then
the spectrum of $H(s,\widetilde{W}(s))$ is the sequence of the eigenvalues $\left(\frac{\pi j}{d(s)} \right)^2 -\widetilde{W}(s)-\frac{1}{4a_0^2},\: j=1,2,\dots$. Consequently, for any $\sigma\ge3/2$, the right-hand side of \eqref{op.val.} can be estimated as
\begin{eqnarray*}
\mathrm{tr}\left(H_0-\frac{1}{4a_0^2}\right)_-^\sigma \le\,
2L_{\sigma,1}^{\mathrm{cl}}\int_{s_0}^\infty\sum_{j=1}^\infty\left(-\left(\frac{\pi j}{d(s)} \right)^2 +\widetilde{W}(s)+\frac{1}{4a_0^2}\right)_+^{\sigma+\frac{1}{2}}\,d s.
\end{eqnarray*}
Combining the above bound with (\ref{est.}) and using the unitary equivalence between $H_{\Omega_1}$ and $\widetilde{H}_{\widetilde{\Omega}_1}$, we get
\begin{eqnarray}\nonumber
\mathrm{tr}\left(H_{\Omega_1}-\frac{1}{4a_0^2}\right)_-^\sigma\le\,
2L_{\sigma,1}^{\mathrm{cl}}\int_{s_0}^\infty\sum_{j=1}^\infty\left(-\left(\frac{\pi j}{d(s)} \right)^2 +\widetilde{W}(s)+\frac{1}{4a_0^2}\right)_+^{\sigma+\frac{1}{2}}\,d s\\\nonumber\le2L_{\sigma,1}^{\mathrm{cl}}\int_{s_0}^\infty\sum_{1\le j\le\frac{\sqrt{\widetilde{W}(s)+(4a_0^2)^{-1}}\,d(s)}{\pi}}\left(\widetilde{W}(s)+\frac{1}{4a_0^2}-\left(\frac{\pi j }{d(s)} \right)^2 \right)_+^{\sigma+\frac{1}{2}}\,d s\\\label{first ineq.}\le\frac{2L_{\sigma,1}^{\mathrm{cl}}}{\pi}\int_{s_0}^\infty \sqrt{\widetilde{W}(s)+\frac{1}{4a_0^2}}\, \left(\widetilde{W}(s)+\frac{1}{4a_0^2}-\left(\frac{\pi }{d(s)} \right)^2 \right)_+^{\sigma+\frac{1}{2}} d(s)\,d s.
\end{eqnarray}
We now move on to the $H_{\Omega_2}$ operator.
\begin{figure}[!t]
\begin{center}
\includegraphics[clip, trim=3.5cm 9.5cm 3.5cm 8cm, width=0.75\textwidth]{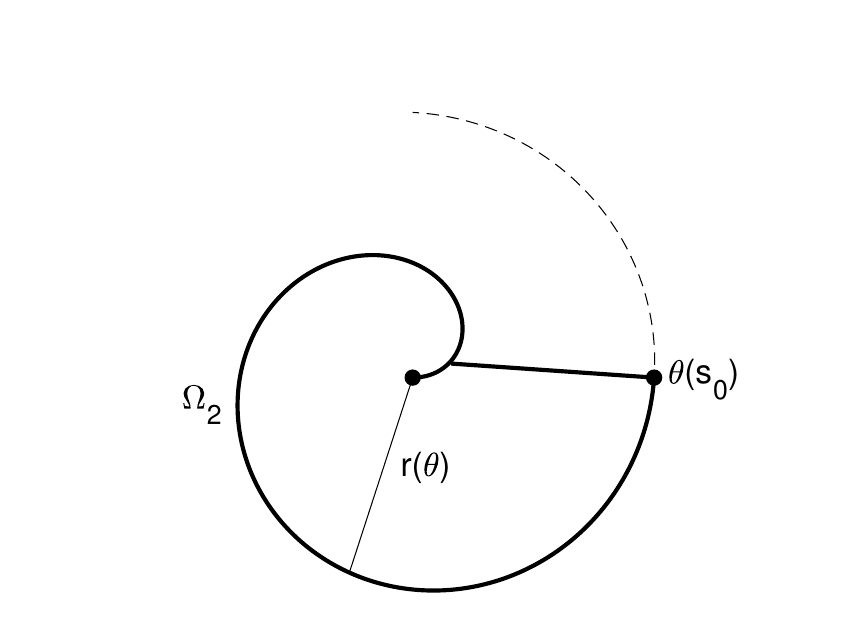}
\end{center}
\caption{The region $\Omega_2$}
\end{figure}
Let $l$ be the segment obtained by extending the interval $\{s=s_0\}\times\{u\in(0, d(s_0))\}$ to the intersection with the boundary of $\Omega_2$ from the left side. It divides the region into two parts; we denote by $\Omega_2^1$ and $\Omega_2^2$ the upper and lower one, respectively. Using the Neumann bracketing we get
\begin{equation}
\label{Omega2}
H_{\Omega_2}\ge H_{\Omega_2^1}\oplus H_{\Omega_2^2},
\end{equation}
 where both operators  $H_{\Omega_2^1}$ and $ H_{\Omega_2^2}$ are the restrictions of $H_{\Omega_2}$ to $\Omega_2^1$ and on $\Omega_2^2$, respectively, with the additional Neumann condition imposed on $l$. Let $\Omega_2^{1,\mathrm{sym}}$ and $\Omega_2^{2,\mathrm{sym}}$ be the unions of $\Omega_2^1$ and $\Omega_2^2$  and their mirror images with respect to $l$, and let $H_{\Omega_2^{1,\mathrm{sym}}}$ and $H_{\Omega_2^{2,\mathrm{sym}}}$ be the Dirichlet Laplacians on them.It is easy to see that one can construct the eigenfunction of $H_{\Omega_2^{1,\mathrm{sym}}}$ and $H_{\Omega_2^{2,\mathrm{sym}}}$ by symmetric continuation of the eigenfunction of $H_{\Omega_2^1}$ and $H_{\Omega_2^2}$ to their mirror images with respect to $l$.
Thus the spectra of $H_{\Omega_2^{1}}$ and $H_{\Omega_2^{2}}$ are contained in the spectra $H_{\Omega_2^{1,\mathrm{sym}}}$ and let $H_{\Omega_2^{2,\mathrm{sym}}}$.
Using the Berezin inequality \cite{Be72a, Be72b} for any $\Lambda>0$ and $\sigma\ge1$ one can estimate the traces of $H_{\Omega_2^{j,\mathrm{sym}}},\: j=1,2,$ as follows
$$
\mathrm{tr}\left(H_{\Omega_2^{j,\mathrm{sym}}}-\Lambda\right)_-^\sigma\le L_{\sigma, 2}^{\mathrm{cl}}\,\mathrm{vol}(\Omega_2^{j,\mathrm{sym}})\, \Lambda^{\sigma+1}, \quad j=1,2,
$$
where $\sigma\ge1$ is a number and the semiclassical constant
\begin{equation}\label{L2}
L_{\sigma, 2}^{\mathrm{cl}}:= \frac{\Gamma(\sigma+1)}{4\pi\Gamma(\sigma+2)} = \frac{1}{4\pi(\sigma+1)}.
\end{equation}

Hence if $\Lambda=\frac{1}{4a_0^2}$ and $\sigma\ge1$ we obtain
\begin{equation}\label{Berezin}
\mathrm{tr}\left(H_{\Omega_2^{j,\mathrm{sym}}}-\frac{1}{4a_0^2}\right)_-^\sigma\le \frac{L_{\sigma, 2}^{\mathrm{cl}}\,\mathrm{vol}(\Omega_2^{j,\mathrm{sym}})}{(4a_0^2)^{\sigma+1}}, \quad j=1,2.
\end{equation}

Finally, combining the inequalities (\ref{bracketing}) and (\ref{first ineq.})--(\ref{Berezin}) we arrive at
\begin{align*}
& \mathrm{tr}\left(H_\Omega-\frac{1}{4a_0^2}\right)_-^\sigma \\[.3em] & \le\,\frac{2L_{\sigma,1}^{\mathrm{cl}}}{\pi}\int_{s_0}^\infty \sqrt{\widetilde{W}(s)+\frac{1}{4a_0^2}}\, \left(\widetilde{W}(s)+\frac{1}{4a_0^2}-\left(\frac{\pi }{d(s)} \right)^2 \right)_+^{\sigma+\frac{1}{2}} d(s)\,d s\\&
 + \frac{L_{\sigma, 2}^{\mathrm{cl}}\,\left(\mathrm{vol}(\Omega_2^{1,\mathrm{sym}})+\mathrm{vol}(\Omega_2^{2,\mathrm{sym}}\right)}{(4a_0^2)^{\sigma+1}}.
\end{align*}
Since
$$
\mathrm{vol}(\Omega_2^{1, \mathrm{sym}})+\mathrm{vol}(\Omega_2^{2, \mathrm{sym}})
=2\mathrm{vol}(\Omega_2^1)+2\mathrm{vol}(\Omega_2^2) =2\mathrm{vol}(\Omega_2)
$$
then the above inequality implies
\begin{align}\nonumber
&\mathrm{tr}\left(H_\Omega-\frac{1}{4a_0^2}\right)_-^\sigma\\\label{second ineq.}& \le\,\frac{2L_{\sigma,1}^{\mathrm{cl}}}{\pi}\int_{s_0}^\infty \sqrt{\widetilde{W}(s)+\frac{1}{4a_0^2}}\, \left(\widetilde{W}(s)+\frac{1}{4a_0^2}-\left(\frac{\pi }{d(s)} \right)^2 \right)_+^{\sigma+\frac{1}{2}} d(s)\,d s
 + \frac{L_{\sigma, 2}^{\mathrm{cl}}\,\mathrm{vol}(\Omega_2)}{4^{\sigma+1/2} a_0^{2\sigma+2}}\end{align}
which completes the proof.
\end{proof}

\begin{remark} \label{remark} 
Note that Theorem \ref{Main} remains valid for smaller powers, $\sigma\ge 1/2$, if we replace the
constants $L_{\sigma,1}^{\mathrm{cl}}$ and $L_{\sigma,2}^{\mathrm{cl}}$ in the right-hand side of \eqref{theorem} by $r(\sigma,1) L_{\sigma,1}^{\mathrm{cl}}$ with the factor $r(\sigma,1)\le2$ if $\sigma< 3/2$, and by $2\left(\frac{\sigma}{\sigma+1}\right)^\sigma\,L_{\sigma,2}^{\mathrm{cl}}$ as follows
\begin{align}\nonumber
&\mathrm{tr}\left(H_\Omega-\Lambda\right)_-^\sigma \\
\nonumber&\le
\,\frac{2r(\sigma,1)  L_{\sigma,1}^{\mathrm{cl}}}{\pi}\int_{s_0}^\infty \sqrt{\widetilde{W}(s)+\frac{1}{4a_0^2}}\, \left(\widetilde{W}(s)+\frac{1}{4a_0^2}-\left(\frac{\pi }{d(s)} \right)^2 \right)_+^{\sigma+\frac{1}{2}} d(s)\,d s
 \\\label{theorem1}&+ \left(\frac{\sigma}{\sigma+1}\right)^\sigma\,L_{\sigma,2}^{\mathrm{cl}}\frac{\,\mathrm{vol}(\Omega_2)}{4^\sigma\,a_0^{2\sigma+2}}.
\end{align} 

The proof follows very closely the proof of Theorem \ref{Main} with the help of (\ref{new result.}) and the Berezin inequality for $0\le\sigma < 1$ \cite{L97}
$$
\mathrm{tr}\left(H_{\Omega_2^{j,\mathrm{sym}}}-\Lambda\right)_-^\sigma\le 2\left(\frac{\sigma}{\sigma+1}\right)^\sigma\,L_{\sigma,2}^{\mathrm{cl}}\,\mathrm{vol}(\Omega_2^{j,\mathrm{sym}})\, \Lambda^{\sigma+1}, \quad j=1,2. 
$$
\end{remark}

\begin{remark}\label{cond.} 
Suppose that functions $\rho'$ and $\rho''$ in (a) are bounded. Let us show that the right-hand sides of the inequalities (\ref{theorem}) and (\ref{theorem1}) are finite for any $\sigma>1/2$, provided by the fact that $(d(s)- 2\pi a_0)_+$ belongs to $L^{\sigma+1/2}(s_0, \infty)$.

Given the formula (\ref{curvature}) it follows that
\begin{equation}\label{gammatheta}
\gamma(\theta)=\frac{1}{a_0 \theta}+\mathcal{O}\left(\frac{1}{\theta^2}\right)\quad\text{ as}\quad \theta\to\infty.
\end{equation}

To express the asymptotic behaviour of $\gamma$ in terms of the variable $s$, we first use the equation
(\ref{s}) to get
\begin{equation}\label{s.}s(\theta)=\frac{a_0}{2}\theta^2+ \mathcal{O}(\theta),\end{equation}
and consequently
$$
\theta=\theta(s)=\sqrt{\frac{2s}{a_0}}+\mathcal{O}(1).
$$

Applying the above expression in (\ref{gammatheta}) we finally get
$$
\gamma(s)=\frac{1}{\sqrt{2a_0 s}}\left(1+\mathcal{O}\left(\frac{1}{\sqrt{s}}\right)\right)\quad\text{ as}\quad s\to\infty.
$$

Since we have the asymptotics for $\gamma(s)$, we find that
\begin{equation}\label{asW}
\widetilde{W}(s)= \frac{1}{8a_0 s}(1+o(1))\quad\text{ as}\quad s\to\infty.
\end{equation}

Given these asymptotics, it is easy to see that the right-hand sides of the inequalities (\ref{theorem}) and (\ref{theorem1}) are finite for any 
$\sigma>1/2$ if
$$
\int_{s_0}^\infty\left(\frac{1}{4a_0^2}-\left(\frac{\pi }{d(s)} \right)^2 \right)_+^{\sigma+\frac{1}{2}}\,d s<\infty,
$$
which, due to the boundedness of $d(s)$, is equivalent to
$$
\int_{s_0}^\infty(d(s)- 2\pi a_0)_+^{\sigma+1/2}\,d s<\infty.
$$ 

\end{remark}

\section{Absence of the infinite discrete spectrum}\label{Absence of the infinite discrete spectrum}

In this section we will describe the class of the asymptotically Archimedean spiral-shaped regions constructed by 
a curve requiring (\ref{assa}) such that the Dirichlet Laplacian $H_\Omega$ has only finite number of the eigenvalues below $\frac{1}{4a_0^2}$. The question of constructing these regions is related to the question of finding the appropriate functions $\rho$ for which we assume that $\rho'$ and $\rho''$ are bounded functions.
We follow the notation of previous chapter and make use of the
estimate (\ref{bracketing}) with some $s_0>0$. Due to the boundedness of the region $\Omega_2$, the spectrum of the operator $H_{\Omega_2}$ is purely discrete and and thus consists of most finite number of eigenvalues below $\frac{1}{4a_0^2}$, then our goal is proved if we establish the absence of the spectrum below $\frac{1}{4a_0^2}$ for the operator $H_{\Omega_1}$. Given (\ref{first ineq.}), this is guaranteed if 
$$
\widetilde{W}(s)+\frac{1}{4a_0^2}\le \left(\frac{\pi}{d(s)}\right)^2
$$
which is equivalent to
$$
\frac{4\pi^2 a_0^2-(d(s))^2}{4a_0^2 (d(s))^2}\ge\widetilde{W}(s)
$$
and finally to
\begin{eqnarray}\nonumber
2\pi a_0-d(s)\ge\frac{\left(2a_0 d(s)\right)^2 \widetilde{W}(s)}{2\pi a_0+d(s)}\\\label{expression}\ge \alpha\widetilde{W}(s),
\end{eqnarray}
where 
$$\alpha=4a_0^2 \left(\inf_{z>s_0}\frac{(d(z))^2}{2\pi a_0+d(z)}\right)>0.$$

Since $\rho'$ and $\rho''$ are bounded functions then (\ref{asW}) is valid. Hence the right hand side of (\ref{expression}) behaves as 
\begin{equation}\label{ass.}
\frac{\alpha}{8a_0 s}(1+\overline{o}(1))\quad\text{ as} \quad s\to\infty.
\end{equation}

On the other hand, using (\ref{s.}) again, one rewrites the asymptotics (\ref{ass.}) in terms of $\theta$ as follows
\begin{equation}
\label{ass}
\frac{\alpha}{4a_0^2 \theta^2}(1+\overline{o}(1))\quad\text{ as} \quad \theta\to\infty.
\end{equation}

Let us now estimate from below the left-hand side of (\ref{expression}) and show that for sufficiently large values of $\theta$ it is substantially larger than (\ref{ass}). Given that 
\begin{eqnarray*}
d(s)=d(s(\theta))\le r(\theta+2\pi)-r(\theta)\\=2\pi a_0+ \rho(\theta+2\pi)-\rho(\theta)\end{eqnarray*}
one gets
\begin{equation}\label{ff}
2\pi a_0-d(s(\theta))\ge -\rho(\theta+2\pi)+\rho(\theta).
\end{equation}

Let us choose $\rho$ so that 
$$\rho'(\theta)=-\frac{c}{\theta^\gamma}+\mathcal{O}\left(\frac{1}{\theta^{\gamma+1}}\right),$$ 
where $c>0$ and $1<\gamma<2$.

Then the left-hand side of (\ref{ff}) for large values of $\theta$ has the following asymptotic behaviour
$$\frac{2\pi c}{\theta^\gamma}+\mathcal{O}\left(\frac{1}{\theta^{\gamma+1}}\right).
$$

Thus, it is easy to see that the above expression is essentially greater than (\ref{ass}) for large values of $\theta$, which proves the validity of the inequality (\ref{expression}).

\section{Local perturbations of the Archimedean spiral}\label{Local perturbations of the Archimedean spiral}

Remark \ref{cond.} establishes the condition for the convergence of the integrals on the right-hand sides of the inequalities (\ref{theorem}) and (\ref{theorem1}) for $\sigma>1/2$. Here we also deal with $\sigma=1/2$.

Suppose that $\mathcal{C}$ coincides with the Archimedean spiral outside of the interval $(s_1, s_2)\subset\mathbb{R}$, where as usual $s$ is the natural parameter of the curve.

As a first step we shall estimate the expression $\frac{1}{4a_0^2}-\left(\frac{\pi }{d(s)} \right)^2$ for $s$ outside the interval $(s_1, s_2)$. Since the spiral $r=r(\theta(s))$ coincides with the Archimedean spiral for such $s$, then using the asymptotic behaviour of the function $\left(\frac{\pi}{d(s)}\right)^2=\left(\frac{\pi}{d(s(\theta))}\right)^2$ with respect to the variable $\theta$, as presented in \cite{ET21}
\begin{align*}&
\left(\frac{\pi}{d(s)}\right)^2= \frac{1}{4a_0^2}+\frac{1}{4a_0^2\theta^2}+\frac{\pi}{2a_0^2\theta^3}+\frac{\pi^2}{a_0^2\theta^4}+\frac{\pi(4\pi^2-1)}{4a_0^2\theta^5}+\mathcal{O}(\theta^{-6})\\&= \frac{1}{4a_0^2}+\frac{1}{4a_0^2\theta^2}+\mathcal{O}(\theta^{-3}) \quad\text{as}\quad \theta\to\infty.
\end{align*}

Hence
\begin{equation}\label{fin.}
\frac{1}{4a_0^2}-\left(\frac{\pi }{d(s)} \right)^2=-\frac{1}{4a_0^2\theta^2}+\mathcal{O}(\theta^{-3}) \quad\text{as}\quad \theta\to\infty.
\end{equation}

Now let us proceed to the study of $\widetilde{W}(s)=\widetilde{W}(s(\theta))$ in terms of the variable $\theta.$ Using the expression (\ref{curvature}) and the fact that the curve $\mathcal{C}$ is a local perturbation of the Archimedean spiral for large values of $\theta$, one obtains
$$\gamma(\theta)=\frac{2+\theta^2}{a_0(1+\theta^2)^{3/2}},$$
which in combination with (\ref{s}) gives
\begin{align*}&
\gamma'(s)=\frac{1}{a_0(1+\theta^2)^{1/2}}\gamma'(\theta),\\&
\gamma''(s)=\frac{1}{a_0^2(1+\theta^2)}\gamma''(\theta)-\frac{\theta}{1+\theta^2} \gamma'(\theta).
\end{align*}

Hence using the expression (\ref{W}) for $\widetilde{W}$ we have
$$
\widetilde{W}(s)=\widetilde{W}(s(\theta))= \frac{1}{4a_0^2 \theta^2}+\mathcal{O}(\theta^{-3}).
$$
Combining the above asymptotics with (\ref{fin.}) we get
$$
\widetilde{W}(s)+\frac{1}{4a_0^2}-\left(\frac{\pi }{d(s)} \right)^2=\widetilde{W}(s(\theta))+\frac{1}{4a_0^2}-\left(\frac{\pi }{d(s(\theta))} \right)^2= \mathcal{O}(\theta^{-3}).
$$

Now let us express $\theta$ in terms of the arc length $s$. Using the fact that the local perturbations of the curve do not affect the asymptotic behaviour of $s$, the equation (\ref{s}) implies
$$
s=\frac{a_0}{2}\theta^2+\mathcal{O}(\theta),
$$
and hence 
$$
\widetilde{W}(s)+\frac{1}{4a_0^2}-\left(\frac{\pi }{d(s)} \right)^2= \mathcal{O}(s^{-3/2}).
$$

This establishes that the integrals on the right-hand side of the inequalities (\ref{theorem}) and (\ref{theorem1}) are finite for any $\sigma\ge1/2$.

\section*{Appendix}

Here we provide a proof of the Theorem \ref{operator valued1} which was skipped in Sec.\ref{Main result}. We use a trick whose idea is credited to Rupert Frank \cite{F22}. Let us introduce the operator 
$$\hat{Q}=-\frac{\partial^2}{\partial\,s^2}\otimes \,I_{L^2(\mathbb{R})}+ H(|s|,\widetilde{W}(|s|))$$
defined on $L^2(\mathbb{R}, L^2(\mathbb{R}))$. It is easy to see that one can construct the eigenfunction of $\hat{Q}$  as a symmetric continuation of the eigenfunction of $Q$ on $L^2(\mathbb{R}, L^2(\mathbb{R}))$. Thus the spectrum of $Q$ is contained in the spectrum of $\hat{Q}$.

The spectrum of the operator $\hat{Q}$ can be estimated by using the operator-valued Lieb-Thirring inequality \cite{LW00} as follows
\begin{equation}\label{LW00}
\mathrm{tr}\Big(\hat{Q}\Big)_-^\sigma\le \,L_{\sigma, 1}^{\mathrm{cl}}\int_{\mathbb{R}}\mathrm{tr}\,(H(|s|, \widetilde{W}(|s|)))_-^{\sigma+1/2}(s)\,d s,\quad \sigma\ge3/2,
\end{equation}
where $L_{\sigma, 1}^{\mathrm{cl}}$ is given in (\ref{LTconstant}).

Hence the above inequality gives
\begin{eqnarray*}
\mathrm{tr}\Big(Q\Big)_-^\sigma\le\mathrm{tr}\Big(\hat{Q}\Big)_-^\sigma\le \,L_{\sigma, 1}^{\mathrm{cl}}\int_0^\infty\mathrm{tr}\,(\hat{H}(|s|, \widetilde{W}(|s|)))_-^{\sigma+1/2}(s)\,d s\\=\,2L_{\sigma, 1}^{\mathrm{cl}}\int_0^\infty\mathrm{tr}\,(H(s, \widetilde{W}(s)))_-^{\sigma+1/2}(s)\,d s\,,\quad \sigma\ge3/2.
\end{eqnarray*}

The second part of the Theorem \ref{operator valued1} is just related to the modification of the Lieb-Thirring inequality for operator-valued potentials, defined on the line $\mathbb{R}$, to the powers
$\sigma\ge 1/2$ and the constant $r(\sigma,1) L_{\sigma,1}^{\mathrm{cl}}$ where $r(\sigma,1)\le 2$ if $\sigma<3/2$ \cite{ELW04}, \textcolor{black}{\cite{HLW00}}, which we used to derive the estimate \eqref{LW00}.


\section*{Statements and Declarations}
\subsection*{Funding} This work is supported in part by Instituto Polit\'ecnico Nacional (grant number SIP20241237).
\bigskip

\subsection*{Competing Interests} The authors have no known competing financial interests or personal relationships that could have appeared to influence the work reported in this paper.
\subsection*{Ethics approval} Not applicable.
\subsection*{Consent to participate} Not applicable. 
\subsection*{Consent for publication} Not applicable.
\subsection*{Author contributions} All authors contributed equally to the manuscript and typed, read, and approved the final form of the manuscript, which is the result of an intensive collaboration. 
\subsection*{Availability of data and material} Not applicable.
\subsection*{Code availability} Not applicable.


\begin{thebibliography}{10}

\bibitem[Be72a]{Be72a}
F.A.~Berezin: Covariant and contravariant symbols of operators, \emph{Izv. Akad. Nauk SSSR Ser. Mat.} \textbf{36} (1972), 1134--167.

\bibitem[Be72b]{Be72b}
F.A.~Berezin: Convex functions of operators, \emph{Mat. Sbornik} (NS) \textbf{36}(130) (1972), 268--276.

\bibitem[BE13]{BE13}
D.~Barseghyan, P.~Exner: Spectral estimates for Dirichlet Laplacians and Schr\"{o}dinger operators on geometrically nontrivial cusps, \emph{J. Spect. Theory} \textbf{3} (4) (2013), 465--484.

\bibitem[BE23]{BE23} D.~ Barseghyan, P.~ Exner:
Spectral estimates for Dirichlet Laplacian on spiral-shaped regions
\emph{J. Spect. Theory} \textbf{13} (2023), 247--261.

\bibitem[BK19]{BK19} D.~Barseghyan, A.~Khrabustovskyi:
Spectral estimates for Dirichlet Laplacian on tubes with exploding twisting velocity, \emph{Operators and Matrices} \textbf{13}(2) (2019), 311--322.

\bibitem[EB14]{EB14}
P.~Exner, D.~Barseghyan: Spectral estimates for Dirichlet Laplacians on perturbed twisted tubes, \emph{Operators and Matrices} \textbf{8}(1) (2014), 167--183.

\bibitem[EK15]{EK15}
P.~Exner, H.~Kova\v{r}\'{\i}k: \emph{Quantum Waveguides}, Springer, Cham 2015.

\bibitem[ELW04]{ELW04}
P.~Exner, H.~Linde, T.~Weidl: Lieb-Thirring inequalities for geometrically induced bound states, \emph{Lett. Math. Phys.} \textbf{70} (2004), 83--95.

\bibitem[ET21]{ET21}
P.~Exner, M.~Tater: Spectral properties of spiral-shaped quantum waveguides, \emph{J. Phys. A: Math. Theor.} \textbf{53} (2020), 5050303

\bibitem[Fr22]{F22} R.L.~Frank, \emph{private communication}

\bibitem[FLW22]{FLW22}
R.L.~Frank, A.~Laptev, T.~Weidl: \emph{Schr\"{o}dinger Operators: Eigenvalues and Lieb-Thirring Inequalities}, Cambridge University Press, to appear

\bibitem[GW11]{GW11}
L.~Geisinger, T.~Weidl: Sharp spectral estimates in domain of infinite volume, \emph{Rev. Math. Phys.} \textbf{23} (2011), 615--641.

\bibitem[HLW00]{HLW00} \textcolor{black}{D.~Hundertmark,  A.~Laptev, T.~Weidl: New bounds on the Lieb-Thirring constants, \emph{Inventiones Mathematicae} \textbf{140} (2000), 693--704.}

\bibitem[JLX15]{JLXZY15}
Jiang Xiao-Jun, Li Xiao-Lin, Xu Xin-Ping, Zhang Hai-Chao, Wang Yu-Zhu: Archimedean-spiral-based microchip ring waveguide for cold atoms, \emph{Chinese Phys. Lett.} \textbf{32} (2015), 020301.

\bibitem[L97]{L97} A.~ Laptev: Dirichlet and Neumann Eigenvalue Problems on Domains in Euclidean Spaces, \emph{J. Func. Anal.} \textbf{151} (1997), 531--545.

\bibitem[LY83]{LY83}
P.~Li, S.T.~Yau: On the Schr\"odinger equation and the eigenvalue problem, \emph{Commun. Math. Phys.} \textbf{88} (1983), 309--318.

\bibitem[Li73]{Li73}
E.H.~Lieb: The classical limit of quantum spin systems, \emph{Commun. Math. Phys.} \textbf{31} (1973), 327--340.

\bibitem[LT76]{LT76}
E.H. Lieb, W.~Thirring: Inequalities for the moments of the eigenvalues of the Schr\"odinger Hamiltonian and their relation to Sobolev inequalities, in \emph{Studies in Math. Phys., Essays in Honor of Valentine Bargmann} (E.~Lieb,
B.~Simon and A.S.~Wightman, eds.); Princeton Univ. Press, Princeton, 1976; pp.~269--330.

\bibitem[LW00]{LW00}
A.~Laptev, T.~Weidl: Sharp Lieb-Thirring inequalities in high dimensions, \emph{Acta Math.} \textbf{184} (2000), 87--100.

\bibitem[Ro72]{Ro72} \textcolor{black}{V. G.~Rozenbljum: The eigenvalues of the first boundary value problem in unbounded domains, \emph{Math. Sbornik} \textbf{18} (1972), 235--248.}

\bibitem[Ro73]{Ro73} \textcolor{black}{V. G.~Rozenbljum: The calculation of the spectral asymptotics for the laplace operator in domains of infinite measure,  \emph{Problems of Mathematical Analysis} \textbf{4} (1973),  95--106.}

\bibitem[Re48]{Re48}
\textcolor{black}{F. Rellich: Das Eigenwert problem von $\Delta\lambda+\lambda\mu=0$ in Halbr\"ohren, in ``Studies and Essays'', Interscience, New York 1948; pp.~329--344}

\bibitem[RS78]{RS78}
M.~Reed, B.~Simon: \emph{Methods of Modern Mathematical Physics, IV.~Analysis of Operators}, Academic Press, New York 1978.

\bibitem[Si83]{Si83}\textcolor{black}{ B.~Simon: Non-classical eigenvalue asymptotics, \emph{J. Funct. Anal.} \textbf{53} (1983), 84--98.}

\bibitem[vdB84]{vdB84} \textcolor{black}{M.~van~den Berg: On the spectrum of the Dirichlet Laplacian for horn-shaped regions in $\mathbb{R}^n$ with
infinite volume,  \emph{J. Funct. Anal.} \textbf{58} (1984), 150--156.}

\bibitem[vdB92]{vdB92} \textcolor{black}{M.~van~den Berg: Dirichlet-Neumann bracketing for horn-shaped regions,  \emph{J. Funct. Anal.} \textbf{104} (1992), 101--120.}

\end{thebibliography}
\end{document}